\documentclass{amsart}
\usepackage{amsfonts}

\usepackage{amsmath}
\usepackage{epsfig}
\usepackage{epsf}

\theoremstyle{plain}

\newtheorem{corollary}{Corollary}[section]

\newtheorem{definition}{Definition}[section]

\newtheorem{lemma}{Lemma}[section]

\newtheorem{remark}{Remark}[section]

\newtheorem{theorem}{Theorem}[section]
\numberwithin{equation}{section}
\input{tcilatex}

\begin{document}
\title[Kneading Theory for Triangular Maps]{Kneading Theory for Triangular Maps}
\author{Diana A. Mendes}
\address{Instituto Superior de Ci\^{e}ncias do Trabalho e da Empresa, Departamento de
M\'{e}todos Quantitativos, Avenida das For\c{c}as Armadas, 1649-026 Lisboa,
Portugal}
\email{diana.mendes@iscte.pt}
\author{J. Sousa Ramos}
\address{Instituto Superior T\'{e}cnico, Departamento de Matem\'{a}tica, Av. Rovisco
Pais 1, 1049-001 Lisboa, Portugal}
\email{sramos@math.ist.utl.pt}

\begin{abstract}
The main purpose of this paper is to present a kneading theory for
two-dimensional triangular maps. This is done by defining a tensor product
between the polynomials and matrices corresponding to the one-dimensional
basis map and fiber map. We also define a Markov partition by rectangles for
the phase space of these maps. A direct consequence of these results is the
rigorous computation of the topological entropy of two-dimensional
triangular maps. The connection between kneading theory and subshifts of
finite type is shown by using a commutative diagram derived from the
homological configurations associated with $m-$modal maps of the interval.
\end{abstract}

\maketitle



\section{Introduction and Preliminary Notions}

\subsection{Introduction}

The techniques of symbolic dynamics have been applied with significant
success in the study of one-dimensional discrete dynamical systems. In
particular, it is well known by now that the kneading sequence (i.e., the
itinerary of the critical point) is a complete topological invariant of a
unimodal map defined on the real unit interval. Over the last few years,
several works have been devoted to generalize these techniques to
two-dimensional dynamics with some success. However the main obstruction to
construct a similar theory for mappings in the real plane is that they lack
critical points in the usual sense and their dynamical space does not seem
to have a natural order which a priori one-dimensional dynamics do have. On
this issue we recall the pruning front conjecture related to the H\'{e}non
and Lozi mappings which can be encountered in the work of Cvitanovic et al. 
\cite{CGP} and Ishii \cite{ISHI}. Their starting point is to regard a
one-folding mapping of the plane as an incomplete horseshoe, and to measure
its incompleteness compared with the full horseshoe in terms of the pruning
front and the primary pruned region which are a two-dimensional analogue of
the kneading sequence.

The main purpose of this paper is to define symbolic dynamics for
two-dimensional triangular maps. Our approach is different from the pruning
techniques, and basically consists of doing some operations (e.g. tensor
products) between the one-dimensional invariants in order to obtain
two-dimensional topological invariants. To the best of our knowledge, we are
not aware of similar works related to kneading theory for triangular maps.
The class of continuous triangular maps we are considering, $T\left(
x,y\right) =\left( f\left( x\right) ,g_{x}\left( y\right) \right) ,$ has the
particularity of admiting finite critical orbits for each one of the
variables $x$, $y$, that is, both maps $f$ and $g_{x}$ are multimodal, which
leads to the result that the transition from dimension one to dimension two
becomes possible. More exactly, we develop a kneading theory and construct a
Markov partition for two-dimensional triangular maps. The kneading
invariants for these maps are obtained in Theorem 2.1 by doing a tensor
product between the kneading matrices and polynomials associated with the
one-dimensional basis map $f$ and fiber map $g$. The Markov partition for
the phase space of the triangular maps (see Theorem 3.1) is simply given by
the cartesian product of the Markov partitions on the $x$ and $y$ lines and
the transition matrix is computed as the Kronecker product of the transition
matrices corresponding to the one-dimensional components of the triangular
map. This is allowed since one of the variables evolves independently of the
other, and in the case of a periodic orbit, the second variable becomes
independent of the first one if the fiber map is considered to be the
composition of the map $g$ at the periodic points of the orbit of map $f$.
Finally, the connection between kneading theory and subshifts of finite type
is shown in Theorem 4.1 by using a commutative diagram derived from the
homological configurations associated with $m-$modal maps of the interval.

This work is based on three important results. The first one is the
structure of periodic orbits of triangular maps due by Kloeden \cite{Kl} and
Alsed\`{a} and Llibre \cite{ALL}.\ The second is based on the notions of
kneading theory for one-dimensional multimodal maps due to Milnor and
Thurston \cite{MT}, and Lampreia and Sousa Ramos \cite{LSR}. The third is
based on the Kronecker product of matrices and on the tensor product of
polynomials over a ring (see, e.g., \cite{Glas}).

An immediate consequence of these results is the exact computation of the
topological entropy \cite{MSR} of two-dimensional triangular maps (see
Corollary 2.1 and 3.1) and through this we can generalize the estimates for
the topological entropy presented in \cite{AKS}. All results are presented
for two-dimensional triangular maps where the basis map is unimodal, the
fiber map is multimodal and the critical points of these maps are finite
periodic points. For the case of eventually periodic and aperiodic orbits
see \cite{MSR2}.

In the scientific literature the triangular maps are also frequently called
by skew-product maps. These maps have applications in geodesic flows on
Riemannian surfaces of constant negative curvature, in the study of strange
attractors and in certain polynomial endomorphisms of $\mathbb{C}^{n}$. Some
mathematical models which provide triangular maps are also found in physics
and economics. The triangular maps (or skew-product maps) were largely
studied in the context of ergodic properties.

This paper is organized as follows: Section 1 reviews some basic facts on
triangular maps and tensor products. The main results concerning kneading
theory and Markov partitions for triangular maps are presented in Section 2
and 3. In Section 4 we present the connection between kneading theory and
subshifts of finite type and finally in Section 5 some examples are provided.

\subsection{Preliminary notions}

Let $X,Y$ be compact intervals of the real line. A two-dimensional
triangular map is a continuous map of the form 
\begin{equation*}
T\left( x,y\right) =\left( f\left( x\right) ,g\left( x,y\right) \right)
=\left( f\left( x\right) ,g_{x}\left( y\right) \right) ,
\end{equation*}
where $T:X\times Y\rightarrow X\times Y$ splits the rectangle $X\times Y$ in
one-dimensional fibers $Y_{x}=Y\left( x\right) =\left\{ x\right\} \times Y$
for any $x\in X$ such that each fiber is mapped by $T$ in a fiber. The map $%
f $ is called the basis map and $g$ is called the fiber map. If we consider
that $X=Y=I$, where $I$ is a compact interval of the real line, then the set
of all continuous maps from $I^{2}$ into itself will be denoted by $%
C_{\Delta }\left( I^{2},I^{2}\right) .$

Let $P=\left\{ x_{0},x_{1},\ldots ,x_{p-1}\right\} $ be a $p-$periodic orbit
of $f$ such that $f\left( x_{i}\right) =x_{i+1},$ for $i=0,\ldots ,p-2$ and $%
f\left( x_{p-1}\right) =x_{0}.$ Define $g_{P}:Y\rightarrow Y$ by 
\begin{equation}
g_{P}\left( y\right) =g\left( x_{p-1},g\left( x_{p-2},\ldots ,g\left(
x_{1},g\left( x_{0},y\right) \right) ...\right) \right) .  \label{EqGP}
\end{equation}
If $Q=\left\{ y_{0},y_{1},\ldots ,y_{q-1}\right\} $ is a $q-$periodic orbit
of $g_{P}$ such that $g_{P}\left( y_{i}\right) =y_{i+1}$, for $i=0,\ldots
,q-2$ and $g_{P}\left( y_{q-1}\right) =y_{0},$ then we define the product of 
$P$ by $Q$, denoted by $P\cdot Q,$ as follows. First we define a sequence of 
$pq$ points in $Y$ by setting 
\begin{equation}
t_{ip+j}=\left\{ 
\begin{array}{ll}
y_{i} & \text{ if }j=0 \\ 
&  \\ 
g\left( x_{j-1},t_{ip+j-1}\right) & \text{ if }j=1,2,\ldots ,p-1
\end{array}
\right.  \label{EqPer1}
\end{equation}
for $i=0,1,\ldots ,q-1.$ Now we define 
\begin{equation}
P\cdot Q=\left\{ \left( x_{j},t_{ip+j}\right) :j=0,1,\ldots ,p-1\text{ and }%
i=0,1,\ldots ,q-1\right\} ,  \label{EqPer2}
\end{equation}
or more explicitly 
\begin{equation*}
\begin{array}{cccc}
\left( x_{0},y_{0}\right) & \left( x_{1},g\left( x_{0},y_{0}\right) \right)
& ... & \left( x_{p-1},g\left( x_{p-2},\ldots ,g\left( x_{1},g\left(
x_{0},y_{0}\right) \right) \ldots \right) \right) \\ 
\left( x_{0},y_{1}\right) & \left( x_{1},g\left( x_{0},y_{1}\right) \right)
& ... & \left( x_{p-1},g\left( x_{p-2},\ldots ,g\left( x_{1},g\left(
x_{0},y_{1}\right) \right) \ldots \right) \right) \\ 
\vdots & \vdots & \ddots & \vdots \\ 
\left( x_{0},y_{q-1}\right) & \left( x_{1},g\left( x_{0},y_{q-1}\right)
\right) & ... & \left( x_{p-1},g\left( x_{p-2},\ldots ,g\left( x_{1},g\left(
x_{0},y_{q-1}\right) \right) \ldots \right) \right) .
\end{array}
\end{equation*}
Note that $P\cdot Q\subset X\times Y$ and has cardinality $pq.$

It was shown by Kloeden \cite{Kl} that the order of the coexisting cycles
for interval maps (Sharkovsky's order $">_{s}"$) remains true for continuous
triangular maps, that is, for every $s\in \mathbb{N}\cup \left\{ 2^{\infty
}\right\} $ there exists $T\in C_{\Delta }\left( X\times Y,X\times Y\right) $
such that $Per\left( T\right) =S\left( s\right) ,$ where $Per\left( T\right) 
$ denotes the set of all periods of $T$ and $S\left( s\right) =\left\{ k\in 
\mathbb{N}:s\geq _{s}k\right\} .$ If $Per\left( T\right) =S\left( s\right) $
we say that the triangular map $T$ has type $s$. It is also know from
Kloeden \cite{Kl} that each periodic orbit of $T$ can be decomposed into a
``product'' of periodic orbits of $f$ and $g_{P}$, that is:

\begin{lemma}
Let $T=\left( f,g\right) :X\times Y\rightarrow X\times Y$ be a continuous
triangular map. Then the following hold
\end{lemma}

\begin{enumerate}
\item  \textit{If }$f$\textit{\ has a periodic orbit }$P$\textit{\ and }$%
g_{P} $\textit{\ has a periodic orbit }$Q$\textit{, then }$P\cdot Q$\textit{%
\ is a periodic orbit of }$T.$

\item  \textit{Conversely, each periodic orbit of }$T$\textit{\ can be
obtained as a product of a periodic orbit }$P$\textit{\ of }$f$\textit{\ by
a periodic orbit of }$g_{P}.$
\end{enumerate}

From these results, Alsed\`{a} and Llibre \cite{ALL} obtained a
characterization of the possible sets of periods of triangular maps which is
presented below. Before enunciating the corollary, we specify that $%
Orb\left( F\right) $ denotes the set of all periodic orbits of a map $F,$
and $\left| R\right| $ denotes the period of an orbit $R.$

\begin{corollary}
Let $T=\left( f,g\right) :X\times Y\rightarrow X\times Y$ be a continuous
triangular map. Then 
\begin{equation*}
Per\left( T\right) =\bigcup_{P\in Orb\left( f\right) ,Q\in Orb\left(
g_{P}\right) }\left| P\right| \cdot \left| Q\right| .
\end{equation*}
\end{corollary}

We will use Bowen's definition of topological entropy \cite{Bow}. For a
continuous triangular map $T$ we set $h_{fib}\left( T\right) =\sup_{x\in
X}h\left( T,Y_{x}\right) $ to be the topological entropy of $T$ on the fiber 
$Y_{x}.$ The Bowen formula for topological entropy bounds is satisfied,
i.e., 
\begin{equation*}
\max \left\{ h\left( f\right) ,h_{fib}\left( T\right) \right\} \leq h\left(
T\right) \leq h\left( f\right) +h_{fib}\left( T\right)
\end{equation*}
where $h\left( f\right) ,h\left( T\right) $ denote the topological entropy
for $f$ and $T$. If all the fiber maps $g$ are monotone, then $h\left(
T,Y_{x}\right) =0$ and $h\left( T\right) =h\left( f\right) $ (see, e.g. \cite
{KolSn})$.$ If the basis map $f$ is simple (has at most type $2^{\infty }$),
then $h\left( T\right) =h_{fib}\left( T\right) .$

We know that if $T\in C_{\Delta }\left( X\times Y,X\times Y\right) $ is of
type greater than $2^{\infty }$ then the topological entropy of $T$ is
positive, if $T$ is of type less than $2^{\infty },$ the topological entropy
is zero and if $T$ is of type $2^{\infty }$ both cases are possible (see 
\cite{KolSh}, \cite{Kol} and \cite{BEL}). Moreover it was found (see for
example \cite{AKS}) that there are triangular maps of type $2^{\infty }$
with $h\left( T\right) =\infty $ and triangular maps of type $2^{\infty }$
with zero topological entropy which are strongly chaotic homeomorphism when
restricted to a minimal set (see \cite{FPS}), properties that are impossible
for continuous maps on the interval. A characterization of the lower bounds
of topological entropy for triangular maps depending on the set of periods
is given in \cite{AKS}.

Some of the two-dimensional triangular maps for which our results apply are
the generalizes Baker's transformation $B_{a,b}:I^{2}\rightarrow I^{2},I=%
\left[ 0,1\right] $, given by 
\begin{equation*}
B_{a,b}\left( x_{n},y_{n}\right) =\left\{ 
\begin{array}{ll}
\left( bx_{n},ay_{n}\right) , & \text{if\ \ }x_{n}\leq \dfrac{1}{b} \\ 
&  \\ 
\left( bx_{n}-1,ay_{n}+\left( 1-a\right) \right) , & \text{if\ \ }x_{n}\geq 
\dfrac{1}{b}
\end{array}
,\right.
\end{equation*}
where $0\leq a\leq 1/2,$ $b\in \mathbb{R}$, the twisted horseshoe map $%
H_{a,b}:I^{2}\rightarrow I^{2},I=\left[ 0,1\right] ,$ 
\begin{equation*}
H_{a,b}\left( x_{n},y_{n}\right) =\left\{ 
\begin{array}{ll}
\left( ax_{n},\dfrac{x_{n}}{a}+\dfrac{y_{n}}{b}+\dfrac{1}{2}\right) , & 
\text{if\ \ }0\leq x_{n}\leq \dfrac{1}{a} \\ 
&  \\ 
\left( a\left( 1-x_{n}\right) ,\dfrac{x_{n}}{a}+\dfrac{y_{n}}{b}+\dfrac{1}{2}%
\right) , & \text{if\ \ }\dfrac{1}{a}\leq x_{n}\leq 1
\end{array}
,\right.
\end{equation*}
with $a,b$ real parameters and some type of Kaplan-Yorke map, which is given
by 
\begin{equation*}
K_{a,b,c}\left( x_{n},y_{n}\right) =\left( \left( ax_{n}+b\right) \func{mod}%
1,-cy_{n}+\cos \left( 2\pi x_{n}\right) \right) ,
\end{equation*}
where $a,b,c$ are real parameters.

Next we are going to define some notions related to tensor products of
matrices and polynomials.

\begin{definition}
Let $A,B$ be two matrices of type $\left( m\times n\right) $ and $\left(
p\times q\right) .$ The tensor product of $A$ and $B$ is a matrix $C$ of
type $\left( mp\times nq\right) ,$ represented by $C=A\otimes B$ and defined
by 
\begin{equation*}
C=A\otimes B=\left[ 
\begin{array}{cccc}
a_{11}B & a_{12}B & \cdots & a_{1n}B \\ 
\vdots & \vdots & \ddots & \vdots \\ 
a_{m1}B & a_{m2}B & \cdots & a_{mn}B
\end{array}
\right] .
\end{equation*}
\end{definition}

The tensor product of matrices is not in general commutative, that is $%
A\otimes B\neq B\otimes A,$ and the following relations hold when dimensions
are appropriate 
\begin{equation}
\begin{array}{l}
\left( A\otimes B\right) \cdot \left( C\otimes D\right) =\left( A\cdot
C\otimes B\cdot D\right) \smallskip \\ 
\left( A\otimes B\right) ^{T}=A^{T}\otimes B^{T}\smallskip \\ 
\det \left( A_{n\times n}\otimes B_{m\times m}\right) =\left( \det \left(
A\right) \right) ^{m}\cdot \left( \det \left( B\right) \right) ^{n}\smallskip
\\ 
trace\left( A\otimes B\right) =trace\left( A\right) \cdot trace\left(
B\right) \\ 
P_{A\otimes B}\left( t\right) =P_{A}\left( t\right) \otimes P_{B}\left(
t\right) ,
\end{array}
\label{EqTens}
\end{equation}
where $P_{C}\left( t\right) $ is the characteristic polynomial associated
with a matrix $C.$

\begin{definition}
Let $\mathbb{R}\left[ x\right] ^{\ast }=\mathbb{R}\left[ x\right] \backslash
\left\{ 0\right\} $ denote the set of nonzero polynomials. Given 
\begin{equation*}
f\left( x\right) =f_{m}x^{m}+\ldots +f_{1}x+f_{0}\text{ and }g\left(
x\right) =g_{n}x^{n}+\ldots +g_{1}x+g_{0}
\end{equation*}
in $\mathbb{R}\left[ x\right] ^{\ast }$ with $f_{m}\neq 0$ and $g_{n}\neq 0,$
let $M$ denote the splitting field of $fg$ and let $\alpha _{1},\ldots
,\alpha _{m}$ be the roots of $f$, and $\beta _{1},\ldots ,\beta _{n}$ the
roots of $g$ in $M.$ The tensor product of the polynomials $f$ and $g$ is
denoted by $f\otimes g\in M\left[ x\right] $ and is defined to be the
following polynomial of degree $mn$%
\begin{equation}
\left( f\otimes g\right) \left( x\right)
=f_{m}^{n}g_{n}^{m}\prod_{i=1}^{m}\prod_{j=1}^{n}\left( x-\alpha _{i}\beta
_{j}\right) .  \label{EqProdTens}
\end{equation}
\end{definition}

If $f$ and $g$ are monic polynomials and $F$ and $G$ denote their respective
companion matrices, then 
\begin{equation*}
f\otimes g=\det \left( xI-F\otimes G\right)
\end{equation*}
provides a practical method to compute $f\otimes g.$ The tensor product of
polynomials is commutative.

\section{Kneading Theory}

In what follows we define the symbolic dynamics necessary to develop the
kneading theory of the two-dimensional triangular map $T=\left( f,g\right) .$
We denote by $T_{P}=\left( f,g_{P}\right) $ the triangular map consisting of
the basis map $f$ and the fiber map $g_{P},$ where the latter is defined as
in relation (\ref{EqGP}). Since we know that all the properties of the map $%
T_{P}=\left( f,g_{P}\right) $ pass on to the original triangular map $%
T=\left( f,g\right) ,$ we are going to study only the former map.\ In this
section we are exclusively concerned with periodic orbits and finite
matrices. We consider that the basis map $f$ is unimodal and the fiber map $%
g_{P}$ is $m-$modal. We define the necessary tools for a general $m-$modal
map $F.$

\begin{figure}[tbp]
\vspace{5mm} \centerline{\epsfig{file=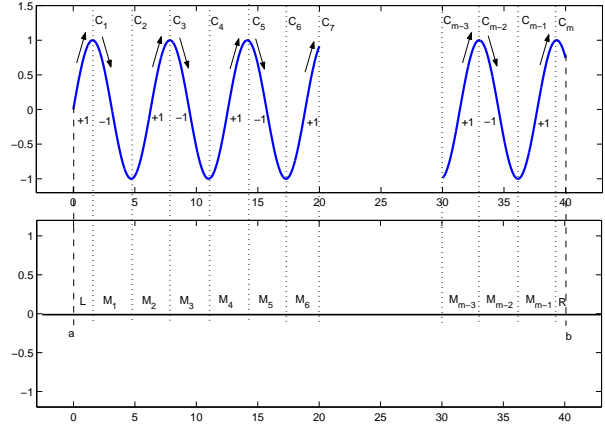,width=8cm}}
\caption{A $m-$modal map.}
\label{FigKNMult}
\end{figure}

Considering that the map $F:X\rightarrow X,$ where $X=\left[ a,b\right] $ is
a compact interval of the real line, is a $m-$modal map and denoting by $%
c_{i},i=1,\ldots ,m$ the $m$ critical points of $F$, we obtain the following
orbits for the critical points for each values of the parameters 
\begin{equation*}
O\left( c_{i}\right) =\left\{ x_{j}^{\left( i\right) }:x_{j}^{\left(
i\right) }=F^{j}\left( c_{i}\right) ,\;j\in \mathbb{N},i=1,\ldots ,m\right\}
.
\end{equation*}
After a reordering of the elements $x_{j}^{\left( i\right) }$ of these
orbits we get a partition $\left\{ X_{k}=\left[ z_{k},z_{k+1}\right]
\right\} $ of the interval $X=\left[ a,b\right] .$ Next, we associate to
each orbit $O\left( c_{i}\right) $ a sequence of symbols $s=S_{1}S_{2}\ldots
S_{j}\ldots ,$ where 
\begin{equation*}
S_{j}=\left\{ 
\begin{array}{ll}
L & \;\text{if }F^{j}\left( c_{i}\right) <c_{1}\medskip \\ 
C_{i} & \;\text{if }F^{j}\left( c_{k}\right) =c_{i},i,k=1,\ldots ,m\medskip
\\ 
M_{i} & \;\text{if }c_{i}<F^{j}\left( c_{k}\right) <c_{i+1},i=1,\ldots
,m-1\medskip \\ 
R & \;\text{if }F^{j}\left( c_{i}\right) >c_{m}.
\end{array}
\right.
\end{equation*}
For simplicity and without lack of generality we assume that the first and
the last critical points of the multimodal map $F$ are maxima. This $m-$%
modal map $F$ and the symbolic partition of the interval $X$ are illustrated
in Figure \ref{FigKNMult}.

The set $\mathcal{A}=\left\{ L,C_{1},M_{1},C_{2},M_{2},\ldots
,C_{m-1},M_{m-1},C_{m},R\right\} $ will be the alphabet of the $m-$modal map 
$F$. We denote the collection of all infinite symbol sequences of $\mathcal{A%
}$ by $\mathcal{A}^{\mathbb{N}}.$ A block or word in $\mathcal{A}$ is a
finite sequence of symbols from $\mathcal{A}$ and the length of a block is
equal to its number of symbols. When $O\left( c_{i}\right) $ is a $k-$%
periodic orbit we get a sequence of symbols that can be characterized by a
block of length $k,$ that is $s^{\left( k\right) }=S_{1}S_{2}\ldots
S_{k-1}C_{i}.$ The shift map $\sigma :\mathcal{A}^{\mathbb{N}}\rightarrow 
\mathcal{A}^{\mathbb{N}}$ is given by $\sigma \left( s\right) =u$ where $%
U_{k}=S_{k+1}.$

In what follows we limit our study to maps for which the orbits of the $m$
critical points are periodic of periods $p_{1},\ldots ,p_{m},$ respectively.
Thus the sequence of symbols corresponding to the itineraries of the $m$
critical points are 
\begin{eqnarray*}
s_{1} &=&S_{11}S_{12}\ldots S_{1p_{1}-1}C_{1}\ldots \medskip \\
s_{2} &=&S_{21}S_{22}\ldots S_{2p_{2}-1}C_{2}\ldots \\
&&\vdots \\
s_{m} &=&S_{m1}S_{m2}\ldots S_{mp_{m}-1}C_{m}\ldots .
\end{eqnarray*}
Let $\left( s_{1}^{\left( p_{1}\right) },\ldots ,s_{m}^{\left( p_{m}\right)
}\right) $ be the $m-$tuple of finite blocks which repeat themselves in the
sequences $s_{1},\ldots ,s_{m}.$ Then the realizable $m-$tuples, as
itineraries of the critical points, are called kneading data and are given
by the following rule \cite{LSR}: if $u^{\left( k\right) }=U_{1}U_{2}\ldots
U_{k}$ is one of the elements of the $m-$tuple with $U_{k}=C_{i},i=1,\ldots
,m$ and $k=p_{i},i=1,\ldots ,m$ then the kneading sequence satisfies the
following points:

\begin{enumerate}
\item  If $U_{i}U_{i+1}\ldots =\sigma ^{i-1}\left( u^{\left( k\right)
}\right) =L$ then $\sigma ^{i}\left( u^{\left( k\right) }\right)
<s_{1}^{\left( p_{1}\right) }$

\item  If $U_{i}U_{i+1}\ldots =\sigma ^{i-1}\left( u^{\left( k\right)
}\right) =M_{j},\;$with $j$ odd, then 
\begin{equation*}
s_{j+1}^{\left( p_{j+1}\right) }<\sigma ^{i}\left( u^{\left( k\right)
}\right) <s_{j}^{\left( p_{j}\right) }
\end{equation*}

\item  If $U_{i}U_{i+1}\ldots =\sigma ^{i-1}\left( u^{\left( k\right)
}\right) =M_{j},\;$with $j$ even, then 
\begin{equation*}
s_{j}^{\left( p_{j}\right) }<\sigma ^{i}\left( u^{\left( k\right) }\right)
<s_{j+1}^{\left( p_{j+1}\right) }
\end{equation*}

\item  If $U_{i}U_{i+1}\ldots =\sigma ^{i-1}\left( u^{\left( k\right)
}\right) =R,$ then $\sigma ^{i}\left( u^{\left( k\right) }\right)
<s_{m}^{\left( p_{m}\right) }.$
\end{enumerate}

Let $V$ be a vector space of dimension $\left( m+1\right) $ defined over the
integers having the formal symbols $\left\{ L,M_{1},\ldots
,M_{m-1},R\right\} $ as a basis, then to each sequence of symbols $%
s=S_{0}S_{1}\ldots S_{j}\ldots $ we can associate a sequence $\theta =\theta
_{0}\theta _{1}\ldots \theta _{j}\ldots $ of vectors from $V$ by setting 
\begin{equation*}
\theta _{j}=\prod_{i=0}^{j-1}\varepsilon \left( S_{i}\right) S_{j},
\end{equation*}
where $j>0$ and 
\begin{equation*}
\theta _{0}=S_{0},\;\varepsilon \left( L\right) =\varepsilon \left(
M_{2k}\right) =1,\;\varepsilon \left( R\right) =\varepsilon \left(
M_{2k+1}\right) =-1,\;\varepsilon \left( C_{i}\right) =0.
\end{equation*}
Choosing a linear order in the vector space $V$ in such a way that the base
vectors satisfy 
\begin{equation*}
L<M_{1}<\ldots <M_{m-1}<R
\end{equation*}
we are able to lexicographically order the sequences $\theta $, that is 
\begin{equation*}
\theta <\theta ^{\prime }\text{ iff }\theta _{0}=\theta _{0}^{\prime
},\ldots ,\theta _{i-1}=\theta _{i-1}^{\prime }\text{ and }\theta
_{i}<\theta _{i}^{\prime }
\end{equation*}
for some integer $i\geq 0.$ Finally, introducing $t$ as an undetermined
variable and taking $\theta _{j}$ as the coefficients of a formal power
series $\theta ,$ we obtain 
\begin{equation*}
\theta =\theta _{0}+\theta _{1}t+\theta _{2}t^{2}+\ldots =\sum_{i=0}^{\infty
}\theta _{i}t^{i}.
\end{equation*}

Milnor and Thurston \cite{MT} also introduced basic invariants called
kneading increments, kneading matrices and kneading determinants. The
kneading increments are formal power series that measure the discontinuity
evaluated at the turning points. For the case of a $m-$modal map we have $m$
kneading increments defined by 
\begin{equation}
\nu _{i}\left( t\right) =\theta _{c_{i}^{+}}\left( t\right) -\theta
_{c_{i}^{-}}\left( t\right) ,\;i=1,\ldots ,m,  \label{EqIK}
\end{equation}
where $\theta $ is the invariant coordinate defined previously\ and 
\begin{equation*}
\theta _{c_{i}^{\pm }}\left( t\right) =\lim_{x\rightarrow c_{i}^{\pm
}}\theta _{x}\left( t\right) .
\end{equation*}

After separating the terms associated with the different symbols in (\ref
{EqIK}) we get 
\begin{equation*}
\nu _{i}\left( t\right) =N_{i1}\left( t\right) L+N_{i2}\left( t\right)
M_{1}+\ldots +N_{im}\left( t\right) M_{m-1}+N_{im+1}\left( t\right) R
\end{equation*}
and from these we can define the kneading matrix by 
\begin{equation*}
N_{F}\left( t\right) =\left[ 
\begin{array}{cccc}
N_{11}\left( t\right) & N_{12}\left( t\right) & \cdots & N_{1m+1}\left(
t\right) \\ 
\vdots & \vdots & \ddots & \vdots \\ 
N_{m1}\left( t\right) & N_{m2}\left( t\right) & \cdots & N_{mm+1}\left(
t\right)
\end{array}
\right] .
\end{equation*}
This is an $m\times \left( m+1\right) $ matrix with entries in the ring $Z%
\left[ \left[ t\right] \right] $ of integer formal power series.$\ $Closely
related is the kneading determinant which is defined from the kneading
matrix according to the following formula 
\begin{equation*}
D_{F}\left( t\right) =\frac{D_{1}\left( t\right) }{1-t}=-\frac{D_{2}\left(
t\right) }{1+t}=\ldots =\left( -1\right) ^{m+1}\frac{D_{m+1}\left( t\right) 
}{1-\left( -1\right) ^{m+1}t},
\end{equation*}
where $D_{i}\left( t\right) $ is the determinant obtained by eliminating the 
$i^{th}$column of the kneading matrix. Finally, we define $d_{F}\left(
t\right) $ by 
\begin{equation}
d_{F}\left( t\right) =D_{F}\left( t\right) \left( 1-t^{p_{1}}\right) \left(
1-t^{p_{2}}\right) \ldots \left( 1-t^{p_{m}}\right) =D_{F}\left( t\right)
P_{cyc}\left( t\right)  \label{EqCyclo}
\end{equation}
where $P_{cyc}\left( t\right) $ is a product of cyclotomic polynomials and $%
p_{1},p_{2},\ldots ,p_{m}$ represent the periodicity of each one of kneading
sequences.

In what follows we suppose that the basis map $f$ is unimodal and the map $%
g_{P}$ is $m-$modal. Let $\mathcal{A}_{x}$ and $\mathcal{A}_{y}$ be the
alphabets corresponding to the basis map and to the fiber map and suppose
that 
\begin{equation*}
s_{x}^{\left( p\right) }\text{ and }\mathbf{u}_{y}^{\left( q\right) }=\left(
u_{1}^{\left( q_{1}\right) },\ldots ,u_{m}^{\left( q_{m}\right) }\right) ,
\end{equation*}
are the associated kneading data. It follows that the triangular map has an
orbit of period $pq$ where $q=q_{1}+\ldots +q_{m}.$

Now we are going to present the main result of this section, that is the
characterization of the kneading invariants for a continuous triangular map.

\begin{theorem}
Let $X,Y$ be compact intervals of the real line and let $T=\left( f,g\right)
:X\times Y\rightarrow X\times Y$ be a continuous triangular map. Assume that
the basis map $f$ has a critical orbit $P$ of period $p$ and the map $g_{P}$
has $m$ critical orbits $Q_{1},\ldots ,Q_{m}$ of periods $q_{1},\ldots
,q_{m} $. Let us denote by $N_{f}\left( t\right) $ and $N_{g_{P}}\left(
t\right) $ the kneading matrices associated with the periodic orbits $P$ and 
$Q_{j},j=1,\ldots ,m.$ Then, the kneading matrix of $T$ is given by the
tensor product of $N_{g_{P}}\left( t\right) $ and $N_{f}\left( t\right) ,$
that is 
\begin{equation*}
\left( N_{T}\left( t\right) \right) _{\left( m\times 2\left( m+1\right)
\right) }=\left( N_{g_{P}}\left( t\right) \right) _{\left( m\times \left(
m+1\right) \right) }\otimes \left( N_{f}\left( t\right) \right) _{\left(
1\times 2\right) },
\end{equation*}
and the kneading-determinant is given by 
\begin{equation*}
D_{T}\left( t\right) =D_{g_{P}}\left( t\right) \otimes D_{f}\left( t\right) .
\end{equation*}
\end{theorem}

\begin{proof}
We will divide the proof in two different cases which are related to the
modality of the fiber map. In the first case the map $g_{P}$ is unimodal and
in the second case it is multimodal.

\textbf{Case 1}: Let $f:X\rightarrow X,g_{P}:Y\rightarrow Y$ be unimodal
maps and let us denote by $c_{x}$ and $c_{y}$ the corresponding critical
points. We also denote by $P$ the critical $p-$period orbit of $f$ and by $Q$
the critical $q-$period orbit of $g_{P}.$

Let $\mathcal{A}_{x}=\left\{ L_{x},C_{x},R_{x}\right\} ,\mathcal{A}%
_{y}=\left\{ L_{y},C_{y},R_{y}\right\} $ be the alphabets corresponding to
the basis map and to the fiber map and suppose that 
\begin{equation*}
s_{x}^{\left( p\right) }=\left( S_{1}\ldots S_{p-1}C_{x}\right) ^{\infty }%
\text{ and }s_{y}^{\left( q\right) }=\left( U_{1}\ldots U_{q-1}C_{y}\right)
^{\infty },
\end{equation*}
with $S_{i}\in \mathcal{A}_{x},i=0,\ldots ,p-1$, $U_{j}\in \mathcal{A}%
_{y},j=0,\ldots ,q-1,$ as the itineraries of $P$ and $Q$, respectively. It
follows that the triangular map has an orbit $P\cdot Q$ of period $pq.$ The
symbolic characterization of this orbit is given by the pair 
\begin{equation*}
\left( \left( s_{x}^{\left( p\right) }\right) ^{\left( q\right) },\left(
s_{y}^{\left( q\right) }\right) ^{\left( p\right) }\right) =\left( \left(
S_{1}\ldots S_{p-1}C_{x}\right) ^{q},\left( U_{1}\ldots U_{q-1}C_{y}\right)
^{p}\right) ^{\infty },
\end{equation*}
and the periodic points are obtained by applying the shift map $\sigma $ to
the pair of symbolic sequences in the following way 
\begin{equation}
\left\{ 
\begin{array}{rr}
\left( \sigma ^{k}\left( s_{x}^{\left( p\right) }\right) ,\sigma
^{k+p}\left( s_{y}^{\left( q\right) }\right) \right) & \text{if }p\neq q%
\text{ and }0\leq k\leq pq-1 \\ 
&  \\ 
\left( \sigma ^{k}\left( s_{x}^{\left( p\right) }\right) ,\sigma
^{k+p+i}\left( s_{y}^{\left( q\right) }\right) \right) & \text{if }p=q\text{
and }0\leq k,i\leq p-1.
\end{array}
\right.  \label{EqSimb}
\end{equation}
We recall that, since $P$ and $Q$ are periodic of periods $p$ and $q$ then
we have that $\sigma ^{k}=\sigma ^{k+p}$ and $\sigma ^{k}=\sigma ^{k+q}$ for 
$k=0,\ldots ,pq-1.$

Since variable $x$ evolves independently of variable $y$ and since the map $%
f $ is unimodal, we can obtain a well defined kneading theory for the basis
map $f$. Denote by $\theta _{x}$ and $\nu _{f}\left( t\right) $ the
invariant coordinate and the kneading increment of the map $f$. Recalling
the relation (\ref{EqGP}) and using the $p$ points of the $P-$orbit of $f$
we transform the fiber map $g$ into an independently one-variable map $%
g_{P}\left( y\right) .$ We have already assumed that the map $g_{P}$ is
unimodal and thus let us denote by $\theta _{y}$ and $\nu _{g_{P}}\left(
t\right) $ the invariant coordinate and the kneading increment of this map.
We recall that all these kneading invariants are formal power series.

Let $Z\left[ \left[ t\right] \right] $ be the ring of formal power series
with integer coefficients, and let $V\left[ \left[ t\right] \right] $ and $W%
\left[ \left[ t\right] \right] $ be the modules consisting of all formal
power series with coefficients in the vector spaces $V$ and $W$. Thus each
series $\theta _{x}=\theta \left( x\right) $ and $\theta _{y}=\theta \left(
y\right) $ is an element of $V\left[ \left[ t\right] \right] $ and $W\left[ %
\left[ t\right] \right] $ respectively, which are free modules with basis $%
\left\{ L_{x},R_{x}\right\} $ and $\left\{ L_{y},R_{y}\right\} .$ In other
words, we can uniquely express $\theta _{x}$ and $\theta _{y}$ as the
following sums 
\begin{equation*}
\theta _{x}=\theta _{1x}L_{x}+\theta _{2x}R_{x}\text{ and }\theta
_{y}=\theta _{1y}L_{y}+\theta _{2y}R_{y},
\end{equation*}
with coefficients $\theta _{1x},\theta _{2x},\theta _{1y},\theta _{2y}$
which are formal power series with integer coefficients.

Then there exist a vector space $V\otimes W$ and a bilinear function $%
\otimes :V\times W\rightarrow V\otimes W$ with the following property: any
bilinear function $B:V\times W\rightarrow U$ to any vector space $U$ over
the integer field can be expressed in terms of $\otimes :V\times
W\rightarrow V\otimes W$ as 
\begin{equation*}
B\left( \theta _{x},\theta _{y}\right) =\left( \theta _{x}\otimes \theta
_{y}\right) T
\end{equation*}
for a unique linear function $T:$ $V\otimes W\rightarrow U.$ So, bilinear
functions from $V\times W$ to $U$ can naturally be considered as
homomorphisms from $V\otimes W$ to $U,$ and conversely (see Lang \cite{L}).

This shows that the invariant coordinate $\theta _{xy}$ of the two
dimensional triangular map is given by 
\begin{eqnarray*}
\theta _{xy} &=&\theta _{x}\otimes \theta _{y}=\left( \theta
_{1x}L_{x}+\theta _{2x}R_{x}\right) \otimes \left( \theta _{1y}L_{y}+\theta
_{2y}R_{y}\right) \\
&=&\theta _{1x}\theta _{1y}L_{x}L_{y}+\theta _{1x}\theta
_{2y}L_{x}R_{y}+\theta _{2x}\theta _{1y}R_{x}L_{y}+\theta _{2x}\theta
_{2y}R_{x}R_{y}.
\end{eqnarray*}
Note that the vector space $U$ has the following base $\left\{
L_{x}L_{y},L_{x}R_{y},R_{x}L_{y},R_{x}R_{y}\right\} $ and the kneading
invariant of the triangular map is now given by 
\begin{equation*}
\nu _{xy}\left( t\right) =N_{11}^{xy}\left( t\right)
L_{x}L_{y}+N_{12}^{xy}\left( t\right) L_{x}R_{y}+N_{13}^{xy}\left( t\right)
R_{x}L_{y}+N_{14}^{xy}\left( t\right) R_{x}R_{y}.
\end{equation*}

It follows immediately that the kneading matrix is of type $\left( 1\times
4\right) $ and is given by the tensor product of the kneading matrices
associated with the one dimensional maps, that is 
\begin{eqnarray*}
\left( N_{T}\left( t\right) \right) _{1\times 4} &=&\left[ 
\begin{array}{r}
N_{11}^{xy}\left( t\right) \medskip \\ 
N_{12}^{xy}\left( t\right) \medskip \\ 
N_{13}^{xy}\left( t\right) \medskip \\ 
N_{14}^{xy}\left( t\right)
\end{array}
\right] ^{T}=\left[ 
\begin{array}{r}
N_{11}^{y}\left( t\right) \otimes N_{11}^{x}\left( t\right) \medskip \\ 
N_{11}^{y}\left( t\right) \otimes N_{12}^{x}\left( t\right) \medskip \\ 
N_{12}^{y}\left( t\right) \otimes N_{11}^{x}\left( t\right) \medskip \\ 
N_{12}^{y}\left( t\right) \otimes N_{12}^{x}\left( t\right)
\end{array}
\right] ^{T}\medskip \\
&=&\left[ 
\begin{array}{rr}
N_{11}^{y}\left( t\right) & N_{12}^{y}\left( t\right)
\end{array}
\right] \otimes \left[ 
\begin{array}{rr}
N_{11}^{x}\left( t\right) & N_{12}^{x}\left( t\right)
\end{array}
\right] \medskip \\
&=&\left( N_{g_{P}}\left( t\right) \right) _{1\times 2}\otimes \left(
N_{f}\left( t\right) \right) _{1\times 2}.
\end{eqnarray*}

The entries of the kneading matrices are polynomials and so, when we
construct the tensor product of these matrices, the new entries of the
resulting matrix are also given by tensor products of polynomials.

The kneading determinant can be computed directly from the kneading matrix $%
N_{T}\left( t\right) $ which is equal to the tensor product of the kneading
determinants of the maps $g_{P}$ and $f,$ that is 
\begin{eqnarray*}
D_{T}\left( t\right) &=&N_{11}^{y}\left( t\right) \otimes N_{11}^{x}\left(
t\right) =N_{11}^{y}\left( t\right) \otimes N_{12}^{x}\left( t\right)
\medskip \\
&=&N_{12}^{y}\left( t\right) \otimes N_{11}^{x}\left( t\right)
=N_{12}^{y}\left( t\right) \otimes N_{12}^{x}\left( t\right)
=D_{g_{P}}\left( t\right) \otimes D_{f}\left( t\right) .
\end{eqnarray*}

\textbf{Case 2}: In this case we suppose that $g_{P}$ is a $m-$modal map.
What is different in this case is just the type of the kneading matrix
associated to the critical orbits of the map $g_{P}$, namely the matrix $%
N_{g_{P}}\left( t\right) $ is of type $\left( m\times \left( m+1\right)
\right) .$ The proof follows in the same way as the first case.
\end{proof}

\begin{corollary}
Let $X,Y$ be compact intervals of the real line and let $T=\left( f,g\right)
:X\times Y\rightarrow X\times Y$ be a continuous triangular map. Then the
topological entropy of $T$ is given by $\log \left( 1/t^{\ast }\right) ,$
where $t^{\ast }$ is the smallest positive solution of the equation $%
D_{T}\left( t\right) =0$, where $D_{T}\left( t\right) $ is the kneading
determinant associated with the kneading matrix $N_{T}\left( t\right) $ from
the previous theorem.
\end{corollary}

\begin{proof}
The proof is immediate by considering the kneading matrix constructed in the
previous theorem, that is 
\begin{equation*}
\left( N_{T}\left( t\right) \right) _{m\times 2\left( m+1\right) }=\left(
N_{g_{P}}\left( t\right) \right) _{m\times \left( m+1\right) }\otimes \left(
N_{f}\left( t\right) \right) _{1\times 2}.
\end{equation*}
It follows that the kneading determinant is given by 
\begin{equation*}
D_{T}\left( t\right) =D_{g_{P}}\left( t\right) \otimes D_{f}\left( t\right) ,
\end{equation*}
and since the topological entropy of the map $f$ is given by $\log \left(
1/t_{x}^{\ast }\right) $ where $t_{x}^{\ast }$ is the smallest positive
solution of $D_{f}\left( t\right) $ and the topological entropy of the map $%
g_{P}$ is given by $\log \left( 1/t_{y}^{\ast }\right) $ where $t_{y}^{\ast
} $ is the smallest positive solution of $D_{g_{P}}\left( t\right) $, it
follows that $\left( 1/t^{\ast }\right) =\left( 1/t_{x}^{\ast }\right)
\left( 1/t_{y}^{\ast }\right) $ is the smallest positive solution of $%
D_{T}\left( t\right) $ (see equation (\ref{EqProdTens})) and $\log \left(
1/t^{\ast }\right) =\log \left( 1/t_{x}^{\ast }\right) +\log \left(
1/t_{y}^{\ast }\right) $ is the topological entropy of the triangular map.
\end{proof}

\begin{remark}
When the orbits of the critical points are eventually periodic, similar
results hold (see \cite{MSR2}). For the aperiodic case we need to use
kneading operators (see also \cite{MSR2}).
\end{remark}

\section{Subshifts of Finite Type}

It is well known that the class of maps studied in the previous section
admits a Markov partition which is determined by the itineraries of the
critical points. Once we have a Markov partition, a subshift of finite type
is determined by a transition matrix. Given a Markov partition $\mathcal{R}%
=\left\{ R_{j}\right\} _{j=1}^{m},$ the transition matrix $A=\left(
a_{ij}\right) $ of type $\left( m\times m\right) $ is defined by 
\begin{equation*}
a_{ij}=\left\{ 
\begin{array}{cc}
1 & \text{if }int\left( f\left( R_{i}\right) \right) \cap int\left(
R_{j}\right) \neq \emptyset \\ 
&  \\ 
0 & \text{if }int\left( f\left( R_{i}\right) \right) \cap int\left(
R_{j}\right) =\emptyset
\end{array}
\right. .
\end{equation*}
The subshift space for $A$ is defined as 
\begin{equation*}
\Sigma _{A}=\left\{ s:\mathbb{N}\rightarrow \left\{ 1,2,\ldots ,m\right\}
:a_{s_{i}s_{i+1}}=1\right\} .
\end{equation*}
Letting $\sigma $ be the shift map on the full $m-$shift, $\Sigma
_{m}=\left\{ 1,2,\ldots ,m\right\} ^{\mathbb{N}},$ define $\sigma
_{A}=\sigma \left| _{\Sigma _{A}}:\Sigma _{A}\rightarrow \Sigma _{A}.\right. 
$

Let us denote by $\mathcal{R}_{x},\mathcal{R}_{y},A_{x},A_{y}$ the Markov
partitions and the transition matrices associated with the one-dimensional
basis map $f$ and fiber map $g_{P}$. Then the following theorem holds:

\begin{theorem}
Let $X,Y$ be compact intervals of the real line and let $T=\left( f,g\right)
\in C_{\Delta }\left( X\times Y,X\times Y\right) $ be a continuous
triangular map. Suppose that the basis map $f$ admits a critical orbit $P$
of \ finite period $p$\ and the fiber map $g_{P}$ admits $m$ critical orbits 
$Q_{1},\ldots ,Q_{m}$ of \ finite periods $q_{1},\ldots ,q_{m}$. Then the
Markov partition of the map $T_{P}$ is given by the cartesian product $%
\mathcal{R}_{x}\times \mathcal{R}_{y}$ and the transition matrix $A$ of $%
T_{P}$ is given by the following tensor product: $A=A_{y}\otimes A_{x}.$
\end{theorem}

\begin{proof}
Suppose that $m=1$ and therefore the maps $f$ and $g_{P}$ are unimodal. Let
us denote by $P$ the period $p$ critical orbit of $f$ and by $Q$ the period $%
q$ critical orbit of $g_{P}.$ Then the interval $X$ (on the $x$ axis of the
real plane) is divided in $p-1$ intervals and the interval $Y$ (on the $y$
axis of the real plane) is divided in $q-1$ intervals. It follows that $%
\mathcal{R}_{x}=\left\{ X_{1},\ldots ,X_{p-1}\right\} $ and $\mathcal{R}%
_{y}=\left\{ Y_{1},\ldots ,Y_{q-1}\right\} $ are the Markov partitions
generated by the two unimodal maps. This way the real plane is divided in $%
\left( p-1\right) \times \left( q-1\right) $ rectangles $R_{i},i=1,\ldots
,\left( p-1\right) \left( q-1\right) $, which define the Markov partition of
the two-dimensional map $T_{P},$ since the image of an initial rectangle $%
R_{i}$ either fully covers a region $R_{j}$ in one iteration or misses it
altogether. This is guarantied by continuity and by the property that
horizontal lines go into horizontal lines and vertical lines go into
vertical lines by following the rules of the maps $f$ and $g_{P}.$

To obtain the transition matrix $A$ of the triangular map we can proceed in
two different ways. The first way is to obtain a $\left( p-1\right) \left(
q-1\right) \times \left( p-1\right) \left( q-1\right) $ matrix directly from
the allowed transitions between the rectangles $R_{i}$ of the Markov
partition, and subsequently one can obtain the same matrix\ by doing the
tensor product of the transition matrices associated with each one of the
one-dimensional maps that compose the triangular map, that is 
\begin{equation*}
\left( A\right) _{\left( p-1\right) \left( q-1\right) \times \left(
p-1\right) \left( q-1\right) }=\left( A_{y}\right) _{\left( q-1\right)
\times \left( q-1\right) }\otimes \left( A_{x}\right) _{\left( p-1\right)
\times \left( p-1\right) }.
\end{equation*}

These follow by the same arguments as those used in the main theorem of the
previous section. The case of multimodal maps follows in the same way.
\end{proof}

\begin{corollary}
Suppose that all the hypothesis of the previous theorem are fulfilled$.$
Then the topological entropy of a continuous triangular map $T=\left(
f,g\right) \in C_{\Delta }\left( X\times Y,X\times Y\right) $ is given by
the sum of the topological entropies of the basis map $f$ and the map $%
g_{P}, $ that is 
\begin{equation*}
h\left( T\right) =h\left( f\right) +h\left( g_{P}\right) .
\end{equation*}
\end{corollary}

\begin{proof}
Since the transition matrix of the map $T_{P}$ is given by the tensor
product of the transition matrices associated with the one-dimensional maps $%
f$ and $g_{P},$ it follows that the characteristic polynomial of the map $%
T_{P}$ is also given by the tensor product of the corresponding
characteristic polynomials of maps $f$ and $g_{P}$ (this follows from
property (\ref{EqTens})\ of the tensor product of matrices). Denoting by $%
\lambda _{x}$ and $\lambda _{y}$ the maximal eigenvalues of the matrices $%
A_{x}$ and $A_{y}$ it follows that 
\begin{equation*}
h\left( T_{P}\right) =h\left( T\right) =\log \left( \lambda _{x}\cdot
\lambda _{y}\right) =\log \lambda _{x}+\log \lambda _{y}=h\left( f\right)
+h\left( g_{P}\right) .
\end{equation*}
This ends the proof.
\end{proof}

\begin{remark}
Between the characteristic polynomial $P_{A}(t)=\det \left( I-tA\right) $ of
the transition matrix $A$ and the kneading determinant of the matrix $%
N_{T}\left( t\right) ,$ $D_{T}\left( t\right) ,$ the following relation
exists 
\begin{equation*}
P_{A}\left( t\right) =D_{T}\left( t\right) P_{cyc}\left( t\right)
=d_{T}\left( t\right) ,
\end{equation*}
where $P_{cyc}\left( t\right) $ is a product of cyclotomic polynomials (see
equation (\ref{EqCyclo})). This confirms the same value for the topological
entropy $h\left( T\right) =\log \left( 1/t^{\ast }\right) $ by using both
methods presented in Corollary 2.1 and Corollary 3.1.
\end{remark}

\section{Connections between Kneading Theory and Subshifts of Finite Type}

In this section we show the connection between kneading theory and subshifts
of finite type by using a commutative diagram derived from the homological
configurations associated with $m-$modal maps of the interval.

For this purpose we consider again the $m-$modal map $F:X\rightarrow X$
which we already presented in Section 2. We assume that $\left(
s_{1}^{\left( p_{1}\right) },\ldots ,s_{m}^{\left( p_{m}\right) }\right) $
is a $m-$modal kneading data, i.e., 
\begin{equation*}
\left( s_{1}^{\left( p_{1}\right) },\ldots ,s_{m}^{\left( p_{m}\right)
}\right) =\left( \left( S_{11}\ldots S_{1p_{1}-1}C_{1}\right) ^{\infty
},\ldots ,\left( S_{m1}\ldots S_{mp_{m}-1}C_{m}\right) ^{\infty }\right) ,
\end{equation*}
and denote by $\mathcal{K}$ the set of all finite $m-$modal kneading data.

Given $\left( s_{1}^{\left( p_{1}\right) },\ldots ,s_{m}^{\left(
p_{m}\right) }\right) \in \mathcal{K}$, let $\left\{ X_{i}\right\}
_{i=1}^{p_{1}+\ldots +p_{m}}$ be the union of the sets $\left\{ \sigma
^{i}\left( s_{1}\right) \right\} _{i=1}^{p_{1}},$ $\ldots $ $,\left\{ \sigma
^{i}\left( s_{m}\right) \right\} _{i=1}^{p_{m}}$ and let $\left\{
x_{i}\right\} _{i=1}^{p_{1}+\ldots +p_{m}}$ denote the points of the
interval with itineraries $\mathcal{I}\left( x_{i}\right) =X_{i}.$ Let $\rho 
$ denote the permutation on $\left\{ 1,2,\ldots ,p_{1}+\ldots +p_{m}\right\} 
$ such that 
\begin{equation*}
a\leq x_{\rho \left( 1\right) }<x_{\rho \left( 2\right) }<\ldots <x_{\rho
\left( p_{1}+\ldots +p_{m}\right) }\leq b
\end{equation*}
and let $z_{i}=x_{\rho \left( i\right) }.$ Finally take the subintervals $%
J_{i}=\left[ z_{i},z_{i+1}\right] ,$ for $i=1,\ldots ,p_{1}+\ldots +p_{m}-1.$
The $m-$modal matrix associated with $\left( s_{1}^{\left( p_{1}\right)
},\ldots ,s_{m}^{\left( p_{m}\right) }\right) $ is a $0,1$-matrix denoted by 
$A$.

Let $\mathcal{C}_{0}$ and $\mathcal{C}_{1}$ be the vector spaces of
dimensions $p_{1}+\ldots +p_{m}$ and $p_{1}+\ldots +p_{m}-1$ of $0-$chains
and $1-$chains spanned by $\left\{ x_{i}\right\} _{i=1}^{p_{1}+\ldots
+p_{m}} $ and $\left\{ J_{i}\right\} _{i=1}^{p_{1}+\ldots +p_{m}-1},$
respectively. Consider the matrix $\varphi $ that maps the basis $\left\{
x_{i}\right\} _{i=1}^{p_{1}+\ldots +p_{m}}$ of $\mathcal{C}_{0}$ onto the
basis $\left\{ J_{i}\right\} _{i=1}^{p_{1}+\ldots +p_{m}-1}$ of $\mathcal{C}%
_{1}$ and take $\eta =\varphi \pi ,$ with $\pi =\pi \left( i,j\right)
=\delta _{\rho \left( i\right) ,j}$ the matrix associated with\ the
permutation $\rho .$ Let $\omega $ be the matrix that represents the
rotation associated with the shift operator in $\mathcal{A}^{\mathbb{N}},$
when restricted to each finite block of the sequence, that is, $\omega
\left( x_{i}\right) =x_{i+1},$ for $i\neq p_{1},\ldots ,p_{1}+\ldots +p_{m},$
and $\omega \left( x_{p_{1}}\right) =x_{1},\omega \left(
x_{p_{1}+p_{2}}\right) =x_{p_{1}+1},\ldots ,\omega \left( x_{p_{1}+\ldots
+p_{m}}\right) =x_{p_{1}+\ldots +p_{m-1}+1}.$ Under these conditions we
obtain an endomorphism $\alpha $ in $\mathcal{C}_{1}$ which is induced from
the commutativity of the following diagram 
\begin{equation*}
\begin{array}{ccccc}
& \mathcal{C}_{0} & \overset{\eta }{\longrightarrow } & \mathcal{C}_{1} & 
\\ 
\omega & \downarrow &  & \downarrow & \alpha . \\ 
& \mathcal{C}_{0} & \underset{\eta }{\longrightarrow } & \mathcal{C}_{1} & 
\end{array}
\end{equation*}
Note that except negative signs the matrix corresponding to the map $\alpha $
is no other than the transition matrix obtained from the admissible
transitions among the subintervals $X_{i}$. Namely, the nonzero elements in
the rows corresponding to the subintervals where the function is decreasing
are equal to $-1$ and in the subintervals where the function is increasing
the nonzero elements are equal to $1.$

We denote by $\beta $ the matrix of order $\left( p_{1}+\ldots
+p_{m}-1\right) \times \left( p_{1}+\ldots +p_{m}-1\right) $ defined by 
\begin{equation*}
\beta =\left[ 
\begin{array}{rrrrr}
I_{n\left( L\right) } & 0 & \ldots & 0 & 0 \\ 
0 & -I_{n\left( M_{1}\right) +1} & \ldots & 0 & 0 \\ 
\vdots & \vdots & \ddots & \vdots & \vdots \\ 
0 & 0 & \ldots & I_{n\left( M_{m-1}\right) +1} & 0 \\ 
0 & 0 & \ldots & 0 & -I_{n\left( R\right) }
\end{array}
\right] ,
\end{equation*}
where $I_{n\left( L\right) },I_{n\left( M_{1}\right) +1},\ldots ,I_{n\left(
M_{m-1}\right) +1},I_{n\left( R\right) }$ are identity matrices of rank $%
n\left( L\right) ,$ $n\left( M_{1}\right) +1,$ $\ldots ,n\left(
M_{m-1}\right) +1$ and $n\left( R\right) ,$ respectively, and where $n\left(
S\right) $ represents de number of symbols $S$ in the sequence.

Calling $\left( u_{1},\ldots ,u_{p_{1}+...+p_{m}}\right) $ the $p_{1}+\ldots
+p_{m}$ sequences of symbols corresponding to the kneading data $\left(
s_{1}^{\left( p_{1}\right) },\ldots ,s_{m}^{\left( p_{m}\right) }\right) $
and calling $U_{i}$ the first symbol of the kneading sequence $u_{i}$, we
associate to this sequence a matrix $\gamma =\gamma \left( i,j\right) $ of
type $\left( p_{1}+\ldots +p_{m}\right) \times \left( p_{1}+\ldots
+p_{m}\right) $ defined as follows:

\begin{enumerate}
\item  We have in the main diagonal 
\begin{equation*}
\gamma \left( i,i\right) =\varepsilon \left( U_{i}\right) ,\;i=2,\ldots
,p_{1}+\ldots +p_{m}
\end{equation*}

\item  In the columns $j=1,p_{1}+1,p_{1}+p_{2}+1,...,p_{1}+...+p_{m-1}+1$ we
have 
\begin{eqnarray*}
\gamma \left( i,j\right) &=&\left\{ 
\begin{array}{cc}
-\varepsilon \left( U_{i}\right) & \text{if }U_{i}<C_{k}=U_{j}\smallskip \\ 
0 & \text{if }U_{i}=C_{k}=U_{j}\smallskip \\ 
\varepsilon \left( U_{i}\right) & \text{if }U_{i}>C_{k}=U_{j}
\end{array}
\right. \\
\text{with }k &=&1,\ldots ,m.
\end{eqnarray*}

\item  All the remaining elements of the matrix are equal to zero.
\end{enumerate}

Let $\Theta =\gamma \omega $ and $A=\beta \alpha ,$ where $\Theta $ is a
matrix of type $\left( p_{1}+\ldots +p_{m}\right) \times \left( p_{1}+\ldots
+p_{m}\right) $ and $A$ is the nonnegative transition matrix of type $\left(
p_{1}+\ldots +p_{m}-1\right) \times \left( p_{1}+\ldots +p_{m}-1\right) .$
Then the following diagram commutes: 
\begin{equation*}
\begin{array}{ccccc}
& \mathcal{C}_{0} & \overset{\eta }{\longrightarrow } & \mathcal{C}_{1} & 
\\ 
\Theta & \downarrow &  & \downarrow & A. \\ 
& \mathcal{C}_{0} & \underset{\eta }{\longrightarrow } & \mathcal{C}_{1} & 
\end{array}
\end{equation*}

We observe that $\eta ^{T}=BiD=\partial ,$ with $B$ the square, integral,
invertible $(p_{1}+\ldots $ $+$ $p_{m})-$ dimensional matrix, given by 
\begin{equation*}
B=\left[ 
\begin{array}{rrrrr}
1 & 0 & \ldots & 0 & 0 \\ 
0 & 1 & \ldots & 0 & 0 \\ 
\vdots & \vdots & \ddots & \vdots & \vdots \\ 
0 & 0 & \ldots & 1 & 0 \\ 
-1 & -1 & \ldots & -1 & 1
\end{array}
\right] ,
\end{equation*}
$i:\mathcal{C}_{1}\hookrightarrow \mathcal{C}_{0}$ the inclusion matrix and $%
D$ the square, integral, invertible (in $\mathbb{Z}$) $(p_{1}+\ldots $ $+$ $%
p_{m}-1)-$ dimensional matrix obtained from $\eta ^{T}$ by removing its $%
(p_{1}+\ldots +$ $p_{m})-$th row. Thus, we have the following commutative
diagram: 
\begin{equation*}
\begin{array}{ccccc}
& \mathcal{C}_{1} & \overset{\partial =\eta ^{T}}{\longrightarrow } & 
\mathcal{C}_{0} &  \\ 
\text{$A^{T}$} & \downarrow &  & \downarrow & \Theta ^{T}. \\ 
& \mathcal{C}_{1} & \underset{\partial =\eta ^{T}}{\longrightarrow } & 
\mathcal{C}_{0} & 
\end{array}
\end{equation*}

Now considering again that the basis map $f$ is a unimodal map with one
critical orbit of period $p$ and $g_{P}$ is a $m-$modal map with $m$
critical orbits of periods $q_{1},\ldots ,q_{m},$ we assume that $\left(
s^{\left( p\right) }\right) \in \mathcal{K}_{x}$ and $\left( u_{1}^{\left(
q_{1}\right) },...,u_{m}^{\left( q_{m}\right) }\right) \in \mathcal{K}_{y}$
are the corresponding kneading data. We denote by $\eta _{i},\partial
_{i},\alpha _{i},\omega _{i},\Theta _{i},A_{i}$ with $i=x,y$ the previously
defined maps and matrices associated with the map $f$ of variable $x$ and
with the map $g_{P}$ of variable $y$.

Since the Markov partition of the triangular map is given by rectangles, we
consider $\mathcal{C}_{1}\times \mathcal{C}_{1}$ to define this partition.
Thus, we have the following commutative diagrams: 
\begin{eqnarray}
&& 
\begin{array}{ccccc}
& \mathcal{C}_{1}\times \mathcal{C}_{1} & \overset{\partial _{y}\otimes
\partial _{x}}{\longrightarrow } & \mathcal{C}_{0}\times \mathcal{C}_{0} & 
\\ 
\left( \alpha _{y}\otimes \alpha _{x}\right) ^{T} & \downarrow &  & 
\downarrow & \left( \omega _{y}\otimes \omega _{x}\right) ^{T} \\ 
& \mathcal{C}_{1}\times \mathcal{C}_{1} & \underset{\partial _{y}\otimes
\partial _{x}}{\longrightarrow } & \mathcal{C}_{0}\times \mathcal{C}_{0} & 
\end{array}
\label{EqDiaCom} \\
&&\text{and}  \notag \\
&& 
\begin{array}{ccccc}
& \mathcal{C}_{1}\times \mathcal{C}_{1} & \overset{\partial _{y}\otimes
\partial _{x}}{\longrightarrow } & \mathcal{C}_{0}\times \mathcal{C}_{0} & 
\\ 
\left( A_{y}\otimes A_{x}\right) ^{T} & \downarrow &  & \downarrow & \left(
\Theta _{y}\otimes \Theta _{x}\right) ^{T}. \\ 
& \mathcal{C}_{1}\times \mathcal{C}_{1} & \underset{\partial _{y}\otimes
\partial _{x}}{\longrightarrow } & \mathcal{C}_{0}\times \mathcal{C}_{0} & 
\end{array}
\notag
\end{eqnarray}
To show that these diagrams are commutative it is sufficient to show that 
\begin{eqnarray*}
\left( \omega _{y}\otimes \omega _{x}\right) ^{T}\cdot \left( \partial
_{y}\otimes \partial _{x}\right) &=&\left( \partial _{y}\otimes \partial
_{x}\right) \cdot \left( \alpha _{y}\otimes \alpha _{x}\right) ^{T}. \\
\left( \Theta _{y}\otimes \Theta _{x}\right) ^{T}\cdot \left( \partial
_{y}\otimes \partial _{x}\right) &=&\left( \partial _{y}\otimes \partial
_{x}\right) \cdot \left( A_{y}\otimes A_{x}\right) ^{T}
\end{eqnarray*}
This is immediate since 
\begin{eqnarray*}
\left( \omega _{y}\otimes \omega _{x}\right) ^{T}\cdot \left( \partial
_{y}\otimes \partial _{x}\right) &=&\omega _{y}^{T}\partial _{y}\otimes
\omega _{x}^{T}\partial _{x} \\
\left( \partial _{y}\otimes \partial _{x}\right) \cdot \left( \alpha
_{y}\otimes \alpha _{x}\right) ^{T} &=&\partial _{y}\alpha _{y}^{T}\otimes
\partial _{x}\alpha _{x}^{T}
\end{eqnarray*}
and since $\omega _{y}^{T}\partial _{y}=\partial _{y}\alpha _{y}^{T}$ and $%
\omega _{x}^{T}\partial _{x}=\partial _{x}\alpha _{x}^{T}$ it follows that
the diagram is commutative. The same reasoning can be applied to the second
commutativity diagram.

Considering all the properties presented above we have the following theorem:

\begin{theorem}
For each kneading data $\left( s^{\left( p\right) }\right) \in \mathcal{K}%
_{x}$ and $\left( u_{1}^{\left( q_{1}\right) },...,u_{m}^{\left(
q_{m}\right) }\right) \in \mathcal{K}_{y}$ corresponding to a continuous
triangular map $T\left( x,y\right) =\left( f\left( x\right) ,g\left(
x,y\right) \right) $ we have that: 
\begin{equation*}
D_{T}\left( t\right) \cdot P_{cyc}\left( t\right) =P_{A}\left( t\right)
=P_{\Theta }\left( t\right) ,
\end{equation*}
where $P_{A}\left( t\right) $ is the characteristic polynomial associated
with the transition matrix $A=A_{y}\otimes A_{x}$ and $P_{\Theta }\left(
t\right) $ is the characteristic polynomial associated with the matrix $%
\Theta =\Theta _{y}\otimes \Theta _{x}$.
\end{theorem}

\begin{proof}
The proof it follows immediately since it is a consequence of the
commutativity of the diagram presented in (\ref{EqDiaCom}) and of the
application to this diagram of some results from homological algebra (see 
\cite{L}).
\end{proof}

\section{Examples}

We consider the following two-parameter continuous triangular map 
\begin{equation*}
T_{a,b}\left( x,y\right) =\left( f\left( x\right) ,g\left( x,y\right)
\right) =\left( 1-ax^{2},x-by^{2}\right) .
\end{equation*}
The basis map $f$ and the fiber map $g$ are both represented by quadratic
maps. We fix $a=1.76$, for which the quadratic basis map $f\left( x\right)
=1-ax^{2}$ has a period $3$ orbit given by 
\begin{equation*}
x_{1}=-0.7589,x_{2}=-0.0135,x_{3}=0.9997.
\end{equation*}

\begin{figure}[tbp]
\vspace{5mm} \centerline{\epsfig{file=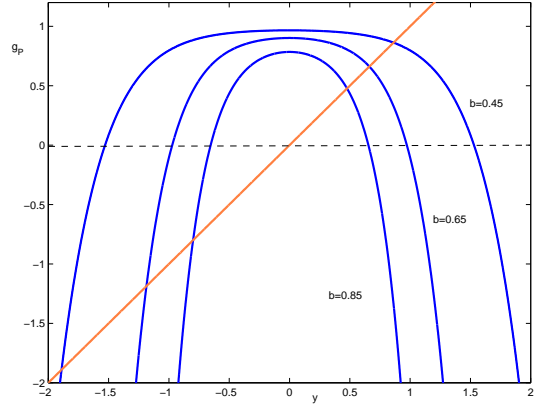,width=7cm}}
\caption{The map $g_{P}$ for $a=1.76$ and different values of $b$.}
\label{FigGP}
\end{figure}

The map $g_{P}\left( y\right) =g\left( x_{3},g\left( x_{2},g\left(
x_{1},y\right) \right) \right) $ has the form 
\begin{equation*}
g_{P}\left( y\right) =0.9997-b(-0.0135-b(-0.7589-by^{2})^{2})^{2}
\end{equation*}
with graphical representation given in Figure \ref{FigGP}\ for different
values of the parameter $b$.

\begin{figure}[tbp]
\vspace{5mm} \centerline{\epsfig{file=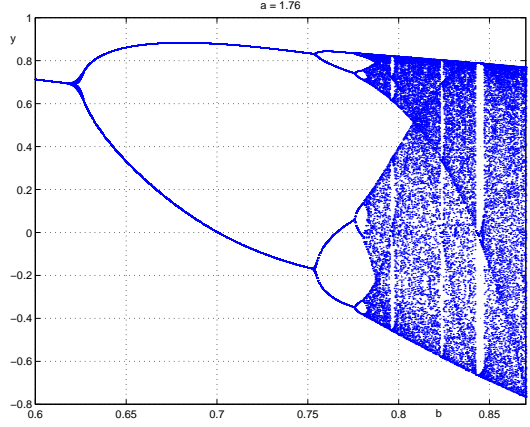,width=7cm}}
\caption{The bifurcation diagram for the map $g_{P}$.}
\label{FigBifGP}
\end{figure}

\begin{figure}[tbp]
\vspace{5mm} \centerline{\epsfig{file=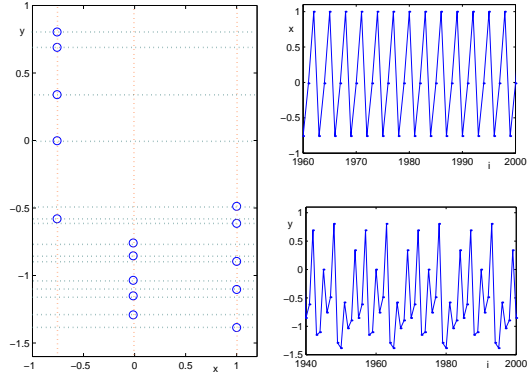,width=7cm}}
\caption{The period $15$ orbit of the triangular map $T$.}
\label{FigP15T}
\end{figure}

\begin{figure}[tbp]
\vspace{5mm} \centerline{\epsfig{file=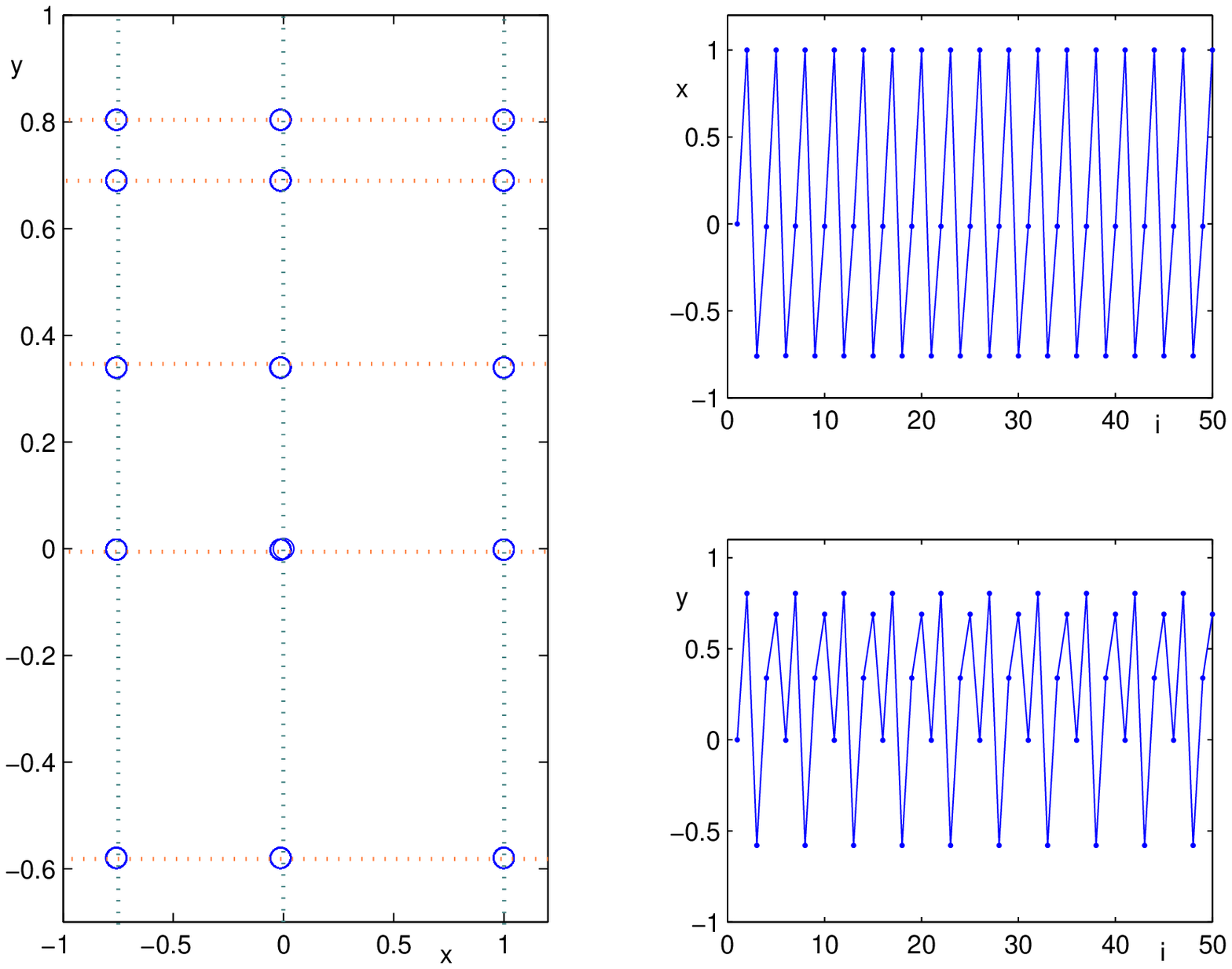,width=7cm}}
\caption{The period $15$ orbit of $\left( f,g_{P}\right) $ for $b=0.823$.}
\label{FigP15GP}
\end{figure}

The bifurcation diagram of the map $g_{P}$ is presented in Figure \ref
{FigBifGP}, when $b$ is varied between $0.6$ and $0.87.$ The map $g_{P}$ has
a unique critical point for $y_{c}=0$. This permits us to define a Markov
partition and a symbolic coding ($L$ if $y<0$ and $R$ if $y>0$). When $%
b=0.823$ the map $g_{P}$ has a period $5$ orbit given by 
\begin{eqnarray*}
y_{1} &=&-0.0018,y_{2}=0.8041,y_{3}=-0.5795, \\
y_{4} &=&0.3396,y_{5}=0.6899,
\end{eqnarray*}
so the maps $T_{P}=\left( f,g_{P}\right) $ and $T=\left( f,g\right) $ have
cycles of period $15$ (Figures \ref{FigP15T}\ and \ref{FigP15GP}).

The kneading sequence of the map $f$ for $a=1.76$ is $\left( RLC\right)
^{\infty }$ which generates a two interval Markov partition on the $x$ line,
with the transition matrix given by 
\begin{equation*}
A_{x}=\left[ 
\begin{array}{cc}
0 & 1 \\ 
1 & 1
\end{array}
\right] _{\left( 2\times 2\right) },
\end{equation*}
and the kneading sequence of the map $g_{P}$ when $b=0.823$ is $\left(
RLRRC\right) ^{\infty }$, which generates a $4$ interval partition on the $y$
line whose transition matrix is given by 
\begin{equation*}
A_{y}=\left[ 
\begin{array}{cccc}
0 & 0 & 1 & 1 \\ 
0 & 0 & 0 & 1 \\ 
0 & 1 & 1 & 0 \\ 
1 & 0 & 0 & 0
\end{array}
\right] _{\left( 4\times 4\right) }.
\end{equation*}
The characteristic polynomials associated with the transition matrices are 
\begin{equation*}
P_{A_{x}}\left( t\right) =1-t-t^{2}\text{ and }P_{A_{y}}\left( t\right)
=1-t-t^{2}+t^{3}-t^{4},
\end{equation*}
and consequently the spectral radius of the matrices $A_{x}$ and $A_{y}$ are
given by $\lambda _{x}=1/t_{x}=\frac{1+\sqrt{5}}{2}\approx 1.6183$ and $%
\lambda _{y}=1/t_{y}\approx 1.5128$ respectively, where $t_{x}$ and $t_{y}$
are the smallest positive solutions of the characteristic polynomials.

\begin{figure}[tbp]
\vspace{5mm} \centerline{\epsfig{file=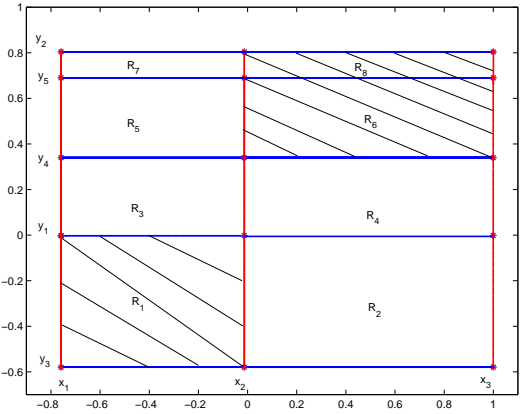,width=7cm}}
\caption{Markov partition for the period 15 orbit of $\left( f,g_{P}\right) $
when $b=0.823$.}
\label{FigRectMar}
\end{figure}

The Markov partition of the map $T_{P}$ is illustrated in Figure \ref
{FigRectMar}. The transitions of the $8$ rectangles $R_{i},i=1,...,8$ that
form the Markov partition of the phase space of map $T_{P}$ are presented by
the following matrix 
\begin{equation*}
A=A_{y}\otimes A_{x}=\left[ 
\begin{array}{cccccccc}
0 & 0 & 0 & 0 & 0 & 1 & 0 & 1 \\ 
0 & 0 & 0 & 0 & 1 & 1 & 1 & 1 \\ 
0 & 0 & 0 & 0 & 0 & 0 & 0 & 1 \\ 
0 & 0 & 0 & 0 & 0 & 0 & 1 & 1 \\ 
0 & 0 & 0 & 1 & 0 & 1 & 0 & 0 \\ 
0 & 0 & 1 & 1 & 1 & 1 & 0 & 0 \\ 
0 & 1 & 0 & 0 & 0 & 0 & 0 & 0 \\ 
1 & 1 & 0 & 0 & 0 & 0 & 0 & 0
\end{array}
\right] _{\left( 8\times 8\right) }.
\end{equation*}
The characteristic polynomial of the transition matrix is given by 
\begin{equation*}
P_{A}\left( t\right) =1-t-4t^{2}+3t^{3}-3t^{4}-5t^{5}+2t^{6}+t^{7}+t^{8},
\end{equation*}
and it follows that the spectral radius of the matrix $A$ is given by $%
\lambda =2.4478=\lambda _{x}\cdot \lambda _{y}$ and consequently the
topological entropy of the maps $T_{P}=\left( f,g_{P}\right) $ and $T=\left(
f,g\right) $ is given by 
\begin{eqnarray*}
h\left( T\right) &=&h\left( T_{P}\right) =\log \left( \lambda \right) =\log
\left( \lambda _{x}\cdot \lambda _{y}\right) \\
&=&\log \left( \lambda _{x}\right) +\log \left( \lambda _{y}\right) =0.8952.
\end{eqnarray*}

The kneading matrix $N_{f}\left( t\right) $ associated with the period $3$
orbit $\left( RLC\right) ^{\infty }$ of the map $f$ is given by 
\begin{equation*}
N_{f}\left( t\right) =\left[ 
\begin{array}{ccc}
-1+\dfrac{2t^{2}}{1-t^{3}} &  & \dfrac{1-2t+t^{3}}{1-t^{3}}
\end{array}
\right] _{\left( 1\times 2\right) }
\end{equation*}
and the kneading matrix $N_{g_{P}}\left( t\right) $ associated with the
period 5 orbit $\left( RLRRC\right) ^{\infty }$ of the map $g_{P}$ is given
by 
\begin{equation*}
N_{g_{P}}\left( t\right) =\left[ 
\begin{array}{ccc}
-1+\dfrac{2t^{2}}{1-t^{5}} &  & \dfrac{1-2t+2t^{3}-2t^{4}+t^{5}}{1-t^{5}}
\end{array}
\right] _{\left( 1\times 2\right) }.
\end{equation*}
Then, the kneading matrix of the triangular map is given by 
\begin{eqnarray*}
N_{T}\left( t\right) &=&N_{g_{P}}\left( t\right) \otimes N_{f}\left( t\right)
\\
&=&\left[ 
\begin{array}{l}
\left( -1+\dfrac{2t^{2}}{1-t^{5}}\right) \otimes \left( -1+\dfrac{2t^{2}}{%
1-t^{3}}\right) \medskip \medskip \\ 
\left( -1+\dfrac{2t^{2}}{1-t^{5}}\right) \otimes \left( \dfrac{1-2t+t^{3}}{%
1-t^{3}}\right) \medskip \medskip \\ 
\left( \dfrac{1-2t+2t^{3}-2t^{4}+t^{5}}{1-t^{5}}\right) \otimes \left( -1+%
\dfrac{2t^{2}}{1-t^{3}}\right) \medskip \medskip \\ 
\left( \dfrac{1-2t+2t^{3}-2t^{4}+t^{5}}{1-t^{5}}\right) \otimes \left( 
\dfrac{1-2t+t^{3}}{1-t^{3}}\right)
\end{array}
\right] ^{T}
\end{eqnarray*}
and the kneading determinant is 
\begin{equation*}
D_{T}\left( t\right) =\left( -1+\dfrac{2t^{2}}{1-t^{5}}\right) \otimes
\left( -1+\dfrac{2t^{2}}{1-t^{3}}\right) ,
\end{equation*}
\begin{equation*}
D_{T}\left( t\right) =\frac{%
1-t-4t^{2}+3t^{3}-3t^{4}-5t^{5}+2t^{6}+t^{7}+t^{8}}{\left( 1-t^{5}\right)
\left( 1-t^{3}\right) }.
\end{equation*}

The smallest positive solution of $D_{T}\left( t\right) =0$ is given by $%
t^{\ast }\approx 0.408515$ and so the topological entropy is computed to be 
\begin{equation*}
h\left( T\right) =\log \left( \frac{1}{t^{\ast }}\right) =\log \left(
2.4478\right) \approx 0.8952.
\end{equation*}
It is easy to observe that 
\begin{equation}
D_{T}\left( t\right) \cdot P_{cyc}\left( t\right) =D_{T}\left( t\right)
\cdot \left( 1-t^{5}\right) \left( 1-t^{3}\right) =P_{A}\left( t\right) .
\label{EqMark}
\end{equation}

In what follows we verify the commutative diagram presented in the previous
section for this map. Thus by simple computations we obtain the following
matrices

\begin{equation*}
\eta ^{T}=\partial =\left[ 
\begin{array}{rrrrrrrr}
1 & -1 & -1 & 0 & 0 & 0 & 0 & 0 \\ 
0 & 1 & 0 & -1 & 0 & 0 & 0 & 0 \\ 
-1 & 0 & 1 & 0 & 0 & 0 & 0 & 0 \\ 
0 & 0 & 0 & 0 & 0 & 0 & 1 & -1 \\ 
0 & 0 & 0 & 0 & 0 & 0 & 0 & 1 \\ 
0 & 0 & 0 & 0 & 0 & 0 & -1 & 0 \\ 
-1 & 1 & 0 & 0 & 0 & 0 & 0 & 0 \\ 
0 & -1 & 0 & 0 & 0 & 0 & 0 & 0 \\ 
1 & 0 & 0 & 0 & 0 & 0 & 0 & 0 \\ 
0 & 0 & 1 & -1 & -1 & 1 & 0 & 0 \\ 
0 & 0 & 0 & 1 & 0 & -1 & 0 & 0 \\ 
0 & 0 & -1 & 0 & 1 & 0 & 0 & 0 \\ 
0 & 0 & 0 & 0 & 1 & -1 & -1 & 1 \\ 
0 & 0 & 0 & 0 & 0 & 1 & 0 & -1 \\ 
0 & 0 & 0 & 0 & -1 & 0 & 1 & 0
\end{array}
\right] _{\left( 15\times 8\right) }
\end{equation*}
\begin{equation*}
\alpha =\alpha _{y}\otimes \alpha _{x}=\left[ 
\begin{array}{rrrrrrrr}
0 & 0 & 0 & 0 & 0 & 1 & 0 & 1 \\ 
0 & 0 & 0 & 0 & -1 & -1 & -1 & -1 \\ 
0 & 0 & 0 & 0 & 0 & 0 & 0 & -1 \\ 
0 & 0 & 0 & 0 & 0 & 0 & 1 & 1 \\ 
0 & 0 & 0 & -1 & 0 & -1 & 0 & 0 \\ 
0 & 0 & 1 & 1 & 1 & 1 & 0 & 0 \\ 
0 & -1 & 0 & 0 & 0 & 0 & 0 & 0 \\ 
1 & 1 & 0 & 0 & 0 & 0 & 0 & 0
\end{array}
\right] _{\left( 8\times 8\right) }
\end{equation*}
\begin{equation*}
\beta =\beta _{y}\otimes \beta _{x}=\left[ 
\begin{array}{rrrrrrrr}
1 & 0 & 0 & 0 & 0 & 0 & 0 & 0 \\ 
0 & -1 & 0 & 0 & 0 & 0 & 0 & 0 \\ 
0 & 0 & -1 & 0 & 0 & 0 & 0 & 0 \\ 
0 & 0 & 0 & 1 & 0 & 0 & 0 & 0 \\ 
0 & 0 & 0 & 0 & -1 & 0 & 0 & 0 \\ 
0 & 0 & 0 & 0 & 0 & 1 & 0 & 0 \\ 
0 & 0 & 0 & 0 & 0 & 0 & -1 & 0 \\ 
0 & 0 & 0 & 0 & 0 & 0 & 0 & 1
\end{array}
\right] _{\left( 8\times 8\right) }
\end{equation*}
\begin{equation*}
\omega =\omega _{y}\otimes \omega _{x}=\left[ 
\begin{array}{rrrrrrrrrrrrrrr}
0 & 0 & 0 & 0 & 1 & 0 & 0 & 0 & 0 & 0 & 0 & 0 & 0 & 0 & 0 \\ 
0 & 0 & 0 & 0 & 0 & 1 & 0 & 0 & 0 & 0 & 0 & 0 & 0 & 0 & 0 \\ 
0 & 0 & 0 & 1 & 0 & 0 & 0 & 0 & 0 & 0 & 0 & 0 & 0 & 0 & 0 \\ 
0 & 0 & 0 & 0 & 0 & 0 & 0 & 1 & 0 & 0 & 0 & 0 & 0 & 0 & 0 \\ 
0 & 0 & 0 & 0 & 0 & 0 & 0 & 0 & 1 & 0 & 0 & 0 & 0 & 0 & 0 \\ 
0 & 0 & 0 & 0 & 0 & 0 & 1 & 0 & 0 & 0 & 0 & 0 & 0 & 0 & 0 \\ 
0 & 0 & 0 & 0 & 0 & 0 & 0 & 0 & 0 & 0 & 1 & 0 & 0 & 0 & 0 \\ 
0 & 0 & 0 & 0 & 0 & 0 & 0 & 0 & 0 & 0 & 0 & 1 & 0 & 0 & 0 \\ 
0 & 0 & 0 & 0 & 0 & 0 & 0 & 0 & 0 & 1 & 0 & 0 & 0 & 0 & 0 \\ 
0 & 0 & 0 & 0 & 0 & 0 & 0 & 0 & 0 & 0 & 0 & 0 & 0 & 1 & 0 \\ 
0 & 0 & 0 & 0 & 0 & 0 & 0 & 0 & 0 & 0 & 0 & 0 & 0 & 0 & 1 \\ 
0 & 0 & 0 & 0 & 0 & 0 & 0 & 0 & 0 & 0 & 0 & 0 & 1 & 0 & 0 \\ 
0 & 1 & 0 & 0 & 0 & 0 & 0 & 0 & 0 & 0 & 0 & 0 & 0 & 0 & 0 \\ 
0 & 0 & 1 & 0 & 0 & 0 & 0 & 0 & 0 & 0 & 0 & 0 & 0 & 0 & 0 \\ 
1 & 0 & 0 & 0 & 0 & 0 & 0 & 0 & 0 & 0 & 0 & 0 & 0 & 0 & 0
\end{array}
\right] _{\left( 15\times 15\right) }
\end{equation*}

\begin{equation*}
\gamma =\left[ 
\begin{array}{rrrrrrrrrrrrrrr}
0 & 0 & 0 & 0 & 0 & 0 & 0 & 0 & 0 & 0 & 0 & 0 & 0 & 0 & 0 \\ 
0 & 0 & 0 & 0 & 0 & 0 & 0 & 0 & 0 & 0 & 0 & 0 & 0 & 0 & 0 \\ 
0 & 0 & 0 & 0 & 0 & 0 & 0 & 0 & 0 & 0 & 0 & 0 & 0 & 0 & 0 \\ 
0 & 0 & 0 & 0 & 0 & 0 & 0 & 0 & 0 & 0 & 0 & 0 & 0 & 0 & 0 \\ 
1 & -1 & 0 & -1 & 1 & 0 & 0 & 0 & 0 & 0 & 0 & 0 & 0 & 0 & 0 \\ 
-1 & 0 & 1 & 1 & 0 & -1 & 0 & 0 & 0 & 0 & 0 & 0 & 0 & 0 & 0 \\ 
0 & 0 & 0 & 0 & 0 & 0 & 0 & 0 & 0 & 0 & 0 & 0 & 0 & 0 & 0 \\ 
-1 & 1 & 0 & 0 & 0 & 0 & 1 & -1 & 0 & 0 & 0 & 0 & 0 & 0 & 0 \\ 
1 & 0 & -1 & 0 & 0 & 0 & -1 & 0 & 1 & 0 & 0 & 0 & 0 & 0 & 0 \\ 
0 & 0 & 0 & 0 & 0 & 0 & 0 & 0 & 0 & 0 & 0 & 0 & 0 & 0 & 0 \\ 
1 & -1 & 0 & 0 & 0 & 0 & 0 & 0 & 0 & -1 & 1 & 0 & 0 & 0 & 0 \\ 
-1 & 0 & 1 & 0 & 0 & 0 & 0 & 0 & 0 & 1 & 0 & -1 & 0 & 0 & 0 \\ 
0 & 0 & 0 & 0 & 0 & 0 & 0 & 0 & 0 & 0 & 0 & 0 & 0 & 0 & 0 \\ 
1 & -1 & 0 & 0 & 0 & 0 & 0 & 0 & 0 & 0 & 0 & 0 & -1 & 1 & 0 \\ 
-1 & 0 & 1 & 0 & 0 & 0 & 0 & 0 & 0 & 0 & 0 & 0 & 1 & 0 & -1
\end{array}
\right] _{\left( 15\times 15\right) }.
\end{equation*}

Now we have that $A^{T}=\left( \beta \,a\right) ^{T}$ and $\Theta
^{T}=\left( \gamma \,\omega \right) ^{T},$ that is

\begin{equation*}
A^{T}=\left[ 
\begin{array}{cccccccc}
0 & 0 & 0 & 0 & 0 & 0 & 0 & 1 \\ 
0 & 0 & 0 & 0 & 0 & 0 & 1 & 1 \\ 
0 & 0 & 0 & 0 & 0 & 1 & 0 & 0 \\ 
0 & 0 & 0 & 0 & 1 & 1 & 0 & 0 \\ 
0 & 1 & 0 & 0 & 0 & 1 & 0 & 0 \\ 
1 & 1 & 0 & 0 & 1 & 1 & 0 & 0 \\ 
0 & 1 & 0 & 1 & 0 & 0 & 0 & 0 \\ 
1 & 1 & 1 & 1 & 0 & 0 & 0 & 0
\end{array}
\right] _{\left( 8\times 8\right) }
\end{equation*}
\begin{equation*}
\Theta ^{T}=\left[ 
\begin{array}{rrrrrrrrrrrrrrr}
0 & 0 & 0 & 0 & 0 & 0 & 0 & 0 & 0 & 0 & 0 & 0 & 0 & 0 & -1 \\ 
0 & 0 & 0 & 0 & 0 & 0 & 0 & 0 & 0 & 0 & 0 & 0 & 0 & -1 & 1 \\ 
0 & 0 & 0 & 0 & 0 & 0 & 0 & 0 & 0 & 0 & 0 & 0 & 0 & 1 & 0 \\ 
0 & 0 & 0 & 0 & 0 & 1 & 0 & 0 & -1 & 0 & 0 & 1 & 0 & 0 & 1 \\ 
0 & 0 & 0 & 0 & 1 & -1 & 0 & -1 & 1 & 0 & 1 & -1 & 0 & 1 & -1 \\ 
0 & 0 & 0 & 0 & -1 & 0 & 0 & 1 & 0 & 0 & -1 & 0 & 0 & -1 & 0 \\ 
0 & 0 & 0 & 0 & 0 & -1 & 0 & 0 & 0 & 0 & 0 & 0 & 0 & 0 & 0 \\ 
0 & 0 & 0 & 0 & -1 & 1 & 0 & 0 & 0 & 0 & 0 & 0 & 0 & 0 & 0 \\ 
0 & 0 & 0 & 0 & 1 & 0 & 0 & 0 & 0 & 0 & 0 & 0 & 0 & 0 & 0 \\ 
0 & 0 & 0 & 0 & 0 & 0 & 0 & 0 & 1 & 0 & 0 & 0 & 0 & 0 & 0 \\ 
0 & 0 & 0 & 0 & 0 & 0 & 0 & 1 & -1 & 0 & 0 & 0 & 0 & 0 & 0 \\ 
0 & 0 & 0 & 0 & 0 & 0 & 0 & -1 & 0 & 0 & 0 & 0 & 0 & 0 & 0 \\ 
0 & 0 & 0 & 0 & 0 & 0 & 0 & 0 & 0 & 0 & 0 & -1 & 0 & 0 & 0 \\ 
0 & 0 & 0 & 0 & 0 & 0 & 0 & 0 & 0 & 0 & -1 & 1 & 0 & 0 & 0 \\ 
0 & 0 & 0 & 0 & 0 & 0 & 0 & 0 & 0 & 0 & 1 & 0 & 0 & 0 & 0
\end{array}
\right] _{\left( 15\times 15\right) }
\end{equation*}
Finally, by the commutativity of this diagram it is easy to see that $\Theta
^{T}\partial =\partial A^{T}.$

The characteristic polynomial of the matrix $\Theta $ is given by 
\begin{equation*}
P_{\Theta }\left( t\right)
=1-t-4t^{2}+3t^{3}-3t^{4}-5t^{5}+2t^{6}+t^{7}+t^{8}
\end{equation*}
and it follows that $P_{\Theta }\left( t\right) =P_{A}\left( t\right) .$
Recalling equation (\ref{EqMark}) we have immediately that 
\begin{equation*}
D_{T}\left( t\right) \cdot P_{cyc}\left( t\right) =P_{A}\left( t\right)
=P_{\Theta }\left( t\right) ,
\end{equation*}
which is the relation presented in Theorem 4.1.

\end{document}